\documentclass[conference]{IEEEtran}
\IEEEoverridecommandlockouts
\usepackage{cite}
\usepackage{amsmath,amssymb,amsfonts}
\usepackage{mathtools}
\mathtoolsset{showonlyrefs}
\usepackage{algorithmic}
\usepackage{graphicx}
\usepackage{textcomp}
\usepackage{xcolor}
\def\BibTeX{{\rm B\kern-.05em{\sc i\kern-.025em b}\kern-.08em
    T\kern-.1667em\lower.7ex\hbox{E}\kern-.125emX}}

\usepackage{hyperref}
\usepackage{orcidlink}
    
\newtheorem{theorem}{Theorem}[section]
\newtheorem{lemma}[theorem]{Lemma}
\newtheorem{proposition}[theorem]{Proposition}
\newtheorem{corollary}[theorem]{Corollary}

\renewcommand{\L}{\Lambda}
\renewcommand{\l}{\lambda}
\newcommand{\R}{\mathbb{R}}

\newcommand{\Z}{\mathbb{Z}}
\newcommand{\N}{\mathbb{N}}
\newcommand{\F}{\mathcal{F}}

\begin{document}
\title{
    Operator-isomorphism pairs and Zak transform methods for the study of Gabor systems
    \thanks{This research was funded in whole or in part by the Austrian Science Fund (FWF) [\href{https://doi.org/10.55776/P33217}{10.55776/P33217}].}
}

\author{
    \IEEEauthorblockN{Markus Faulhuber}
    \IEEEauthorblockA{
        \textit{Faculty of Mathematics} \\
        \textit{University of Vienna}\\
        Vienna, Austria\\
        \href{mailto:markus.faulhuber@univie.ac.at}{markus.faulhuber@univie.ac.at} \quad 0000-0002-7576-5724 \orcidlink{0000-0002-7576-5724}
    }
}

\maketitle

\begin{abstract}
   We collect and summarize known results on the unitary equivalence of Gabor systems by pairs of unitary operators and isomorphisms of the time-frequency plane. The methods are then used to study Gabor systems with Hermite functions. We provide new proofs of some known results and we given an outlook on Gabor frames with double over-sampling.
\end{abstract}


\section{Introduction and Notation}
Gabor systems have been introduced by D.~Gabor in \cite{Gab46} and are frequently used for the stable expansion of functions. By stable, we mean that \eqref{eq:frame} holds true. The basic objects defining a Gabor system are the unitary operators of translation (time-shift) and modulation (frequency-shift)
\begin{equation}
    \mathcal{T}_x f(t) = f(t-x)
    \quad \text{ and } \quad
    \mathcal{M}_\omega f(t) = e^{2 \pi i \omega t}.
\end{equation}
They satisfy the following commutation relation
\begin{equation}\label{eq:CR}
    \mathcal{T}_x \mathcal{M}_\omega = e^{-2 \pi i x \omega} \mathcal{M}_\omega \mathcal{T}_x.
\end{equation}
The combination $\pi(z) = \mathcal{M}_\omega \mathcal{T}_x$ is a time-frequency shift by $z = (x,\omega)$. A Gabor system is a function system of the form
\begin{equation}
    \mathcal{G}(g,\L) = \{ \pi(\l) g \mid \l \in \L \subset \R^2 \},
\end{equation}
where $g \in L^2(\R)$ is sufficiently nice and $\L$ a discrete subset of the time-frequency plane $\R^2$. Often $\L$ is assumed to be a (rectangular) lattice. Here, we will study discrete subsets of $\R^2$ which are the union of shifted copies of the integer lattice $\Z^2$. In particular, we consider the case of integer over-sampling of Gabor systems with Hermite functions.

The Fourier transform of a suitable function $f$ is given by
\begin{equation}
    \F f (y) = \widehat{f}(y) = \int_\R f(t) e^{-2 \pi i y t} \, dt.
\end{equation}
We note that $\F$ is unitary on $L^2(\R)$, i.e., we have
\begin{equation}
    \langle \F f, \F g \rangle = \langle f, g \rangle,
    \quad \text{ where }
    \langle f , g \rangle = \int_\R f(t) \overline{g(t)} \, dt
\end{equation}
denotes the inner product on $L^2(\R)$. The $n$-th order Hermite function is given by (cf.~\cite[Chap.~1, \S~7]{Fol89})
\begin{equation}
    h_n(t) = (-1)^n C_n e^{\pi t^2} \frac{d}{dt} e^{-2 \pi t^2}, \quad t \in \R,
\end{equation}
where $C_n = \frac{2^{1/4}}{\sqrt{(2 \pi)^n 2^n n!}}$ is a normalizing constant.

\section{Intertwining of time-frequency shifts}
By the intertwining properties of time-frequency shifts with certain unitary operators we transfer our study to the case of multi-window Gabor systems, which are systems of the form
\begin{equation}
    \mathfrak{G}_M(g_m,\L) = \bigcup_{m=1}^M \mathcal{G}(g_m, \L), \quad M \in \N.
\end{equation}
The associated frame operator has the form \cite[Chap.~8.3]{Gro01}
\begin{equation}
    S_{\mathfrak{G}_M} = \sum_{m=1}^M \sum_{\l \in \L} \langle f, \pi(\l) g_m \rangle \pi(\l) g_m.
\end{equation}
The Gabor system $\mathfrak{G}_M$ is a frame if and only if there exist constants $0 < A \leq B < \infty$ such that for all $f \in L^2(\R)$
\begin{equation}\label{eq:frame}
    A \lVert f \rVert_2^2 \leq \sum_{m=1}^M \sum_{\l \in \L} | \langle f, \pi(\l) g_m \rangle|^2 \leq B \lVert f \rVert_2^2.
\end{equation}
A Gabor system with a finite upper frame bound $B$ is called a Bessel system. Under mild assumptions on $g$, this condition is met and all series converge unconditionally. In the sequel, we will always assume that we have a Bessel system. The frame property is thus satisfied if we can show that $A > 0$.

Now, consider a unitary operator $\mathcal{U}$ on $L^2(\R)$ and an isomorphism $U$ on $\R^2$, which satisfy
\begin{equation}\label{eq:UTF}
    \mathcal{U} \, \pi(z) \, \mathcal{U}^{-1} = c_U(z) \, \pi(U z), \quad |c_U(z)| = 1, \ z \in \R^2.
\end{equation}
We write $\mathfrak{G}_M (g_m,\L) \cong \mathfrak{G}_M (\mathcal{U}g, U \L)$ and say that these Gabor systems are unitarily equivalent in the sense that
\begin{align*}
    & \ \mathcal{U} \, S_{\mathfrak{G}_M(g_m,\L)} \, \mathcal{U}^{-1} f
    = \sum_{m=1}^M \sum_{\l \in \L} \langle \mathcal{U}^{-1} f, \pi(\l) g_m \rangle \, \mathcal{U} \pi(\l) g_m \\
    = & \ \sum_{m=1}^M \sum_{\l \in \L} \langle f, \mathcal{U} \pi(\l) \mathcal{U}^{-1} \mathcal{U} g_m \rangle \, \mathcal{U} \pi(\l) \mathcal{U}^{-1} \mathcal{U} g_m\\
    = & \ \sum_{m=1}^M \sum_{\l \in \L} \langle f, \pi(U \l) \mathcal{U} g_m \rangle \, \pi(U \l) \mathcal{U} g_m
    = S_{\mathfrak{G}_M(\mathcal{U}g_m, U \L)} f
\end{align*}
In particular, they possess the same frame bounds. Note that in the above calculations the phase factor $c_U$, appearing in \eqref{eq:UTF} also appears as complex conjugate $\overline{c_U}$ and, so, cancels.

\section{Valid Operator-isomorphism pairs}
We will now list examples of pairs of unitary operators and bijective mappings which satisfy \eqref{eq:UTF}, which may be found at various places in the literature. According references will be mentioned on the spot. We limit our attention to classical Gabor systems (so $M=1$) for the moment. Since we have
\begin{equation}
    S_{\mathfrak{G}_M} = \sum_{m=1}^M S_{g_m,\L}
\end{equation}
the results easily carry over to the multi-window case.

We start with metaplectic operators and their projections onto the group of symplectic matrices (see, e.g., \cite[Chap.~1.1]{CorRod20}, \cite[Chap.~4]{Fol89}, \cite[Chap.~7]{Gos11}, \cite[Chap.~9.4]{Gro01}). We have the unitary equivalence $\mathcal{U} S_{g,\L} \mathcal{U}^{-1} = S_{\mathcal{U}g, U\L}$ for pairs
\begin{equation}
    (\mathcal{U},U) \in \{ (\mathcal{D}_a, D_a), (\mathcal{V}_q, V_q)\},
\end{equation}
where $\mathcal{D}_a$ and $\mathcal{V}_q$ are a dilation and a chirp, respectively. For $a>0$ and $q \in \R$, these are given by
\begin{align}
    \mathcal{D}_a f(t) = \tfrac{1}{\sqrt{a}} f \left( \tfrac{t}{a} \right),
    \quad
    \mathcal{V}_q f(t) = e^{\pi i q t^2} f(t).
\end{align}
The corresponding isomorphisms are the symplectic matrices
\begin{align}
    D_a = \begin{pmatrix} a & 0 \\ 0 & a^{-1} \end{pmatrix},
    \quad
    V_q = \begin{pmatrix} 1 & 0 \\ q & 1 \end{pmatrix}.
\end{align}

In time-frequency analysis we seek to view a function in dependence of the time- and the frequency-variable. Thus, similar to the idea of the Radon transform, we may think of $f$ and $\widehat{f}$ as a restriction of a function $F$ of two variables, i.e.,
\begin{equation}
    f(t) = F(x,\omega) \big|_{(x,\omega)=(t,0)}
    \ \text{ and } \
    \widehat{f}(y) = F(x,\omega) \big|_{(x,\omega)=(0,y)}
\end{equation}
So, we look at a function defined on the plane only along a horizontal or a vertical line. However, nothing prevents us from also looking at the function $F(x,\omega)$ along a line at a certain angle in the plane. This leads to the fractional Fourier transform (FrFT) (see \cite{Alm94}, \cite[\S~2.3]{Leest01}, \cite{Nam80}), which is given~by
\begin{equation}\label{eq:FrFT}
    \mathcal{F}_r f(s) = \int_\R f(t) k_r(s,t) \, dt,
\end{equation}
where, for $r\in \R$, the integral kernel is
\begin{equation}\label{eq:kr}
    k_r(s,t) = \sqrt{1 - i \cot(r)} \, e^{\pi i \left(\cot(r) s^2 - 2 \csc(r) s t + \cot(r) t^2\right)}.
\end{equation}
The expression is defined whenever $r \notin \pi\Z$. The cases when $r$ is an integer multiple of $\pi$ can also be handled and we have
\begin{align}
    \F_0 f(s) & = f(s),
    & \F_{\frac{\pi}{2}} f(s) & = \widehat{f}(s),\\
    \F_{\pi} f(s) & = f(-s),
    & \F_{\frac{3\pi}{2}} f(s) & = \widehat{f}(-s),
\end{align}
\begin{equation}
    \F_{r_1} \F_{r_2} f(s) = \F_{r_1+r_2} f(s).
\end{equation}
Besides these basic properties, we also know how the FrFT intertwines translations and modulations. We have
\begin{align}
    \F_r \mathcal{T}_x & = e^{\pi i x^2 \sin(r)\cos(r)} \mathcal{M}_{-x \sin(r)} \mathcal{T}_{x \cos(r)} \F_r \label{eq:TxFrFT}\\
    \F_r \mathcal{M}_\omega & = e^{-\pi i \omega^2 \sin(r) \cos(r)} \mathcal{M}_{\omega \cos(r)} \mathcal{T}_{\omega \sin(r)} \F_r. \label{eq:MwFrFT}
\end{align}
In particular, for a time-frequency shift we have
\begin{equation}\label{eq:TFshift_FrFT}
    \F_r \pi(z) \F_r^{-1} = c_r(z) \pi(R_r z),
\end{equation}
where $c_r(z)$ is a phase factor, so $|c_r(z)|=1$, with phase precisely determined by the relations \eqref{eq:CR}, \eqref{eq:TxFrFT}, \eqref{eq:MwFrFT} and
\begin{equation}
    R_r = \begin{pmatrix}
        \cos(r) & \sin(r)\\
        -\sin(r) & \cos(r)
    \end{pmatrix},
\end{equation}
so $R_r \in SO(2,\R)$ is a rotation by $r \in \R$. It follows from relation \eqref{eq:TFshift_FrFT} that the frame operators $S_{\F_r g,\mathcal{R}_r \L}$ and $S_{g,\L}$ are unitarily equivalent and intertwined by the FrFT, i.e.,
\begin{equation}
    \F_r S_{g,\L} \F_r^{-1} = S_{\F_r g, R_r \L}.
\end{equation}
So, we identified the following operator-isomorphism pairs
\begin{equation}
    (\mathcal{D}_a, D_a), \
    (\mathcal{V}_q, V_q), \
    (\mathcal{F}_r, R_r),
\end{equation}
which yield unitarily equivalent Gabor systems. Clearly, the composition of operators and corresponding isomorphisms is also a valid choice. Moreover, any lattice $\L$ of co-volume 1 in $\R^2$ is the image of $\Z^2$ under a composition of $D_a$, $V_q$, and~$R_r$.
\begin{equation}
    \L = R_r V_q D_a \Z^2.
\end{equation}

The case of a shift is special. Adding a phase factor to time-frequency shifts makes them a non-commutative group with $c_1 \pi (z_1) c_2 \pi(z_2) = \widetilde{c} \, \pi(z_1+z_2)$, where $\widetilde{c}$ is determined by \eqref{eq:CR} and $c_1$ and $c_2$. A straight forward computation shows that
\begin{equation}
    S_{\pi(z)g, \L} = \sum_{\l \in \L} \langle f, \pi(\l) \pi(z) g \rangle \, \pi(\l)\pi(z) g = S_{g,\L+z}
\end{equation}
This fact lies at the heart of the investigation of integer over-sampled Gabor systems with the Zak transform.

\section{The Zak transform}
The Zak transform has become a standard tool to study the frame property of Gabor systems. For a function $f$, its Zak transform is given by
\begin{equation}\label{eq:Zak}
    \mathcal{Z}f(x,\omega) = \sum_{k \in \Z} f(k-x) e^{2 \pi i \omega k} = \sum_{k \in \Z} M_\omega T_{x} f(k)
\end{equation}
The Zak transform is useful in the study of Gabor systems due to the following relations \cite[Thm.~8.3.1, Cor.~8.3.2]{Gro01}:
\begin{equation}
    \mathcal{Z} (S_{g,\Z^2} f) = |\mathcal{Z} g|^2 \mathcal{Z} f,
    \quad
    A \leq |\mathcal{Z} g(x,\omega)|^2 \leq B.
\end{equation}
Thus, a multi-window Gabor system over the integer lattice $\Z^2$ is a frame if and only if
\begin{equation}
    A \leq \sum_{m=1}^M |\mathcal{Z} g_m(x,\omega)|^2 \leq B.
\end{equation}
If we seek to study properties of the system $\bigcup_{m=1}^M \mathcal{G}(g_m,\L)$ with lattice $\L = U \Z^2$, then we may simply apply the appropriate operator-isomorphism pair $(\mathcal{U}^{-1},U^{-1})$ and the Zak transform methods become immediately available. Alternatively, we could define a ``lattice Zak transform'' with additional parameters. Sometimes, for example, the Zak transform
\begin{equation}
    \mathcal{Z}_a f(x,\omega) = \sqrt{a} \sum_{k \in \Z} f(a k - x) e^{2 \pi i a k \omega}, \quad a > 0,
\end{equation}
is used to study Gabor systems over rectangular lattices.

\section{Properties of the Zak transform}
For us, an appropriate function space is $W_0(\R)$, the Wiener space of continuous functions (see, e.g., \cite[Chap.~6.1]{Gro01}):
\begin{equation}
    f \in W_0(\R) \subset C \cap L^\infty (\R)
    \ \Longleftrightarrow \
    \sum_{k \in \Z} \lVert T_k \raisebox{2pt}{$\chi$}_{[0,1]} f \rVert_\infty < \infty.
\end{equation}
The Zak transform is quasi-periodic in the sense that
\begin{align}
    \mathcal{Z}f(x+1,\omega) &  = e^{-2 \pi i \omega} \mathcal{Z}f(x,\omega)\\
    \mathcal{Z}f(x,\omega+1) & = \mathcal{Z} f(x,\omega).
\end{align}
Moreover, from \eqref{eq:CR} and \eqref{eq:Zak}, it is easy to observe that
\begin{equation}
    \mathcal{Z}(\mathcal{M}_\eta \mathcal{T}_\xi f)(x,\omega) = e^{-2 \pi i x \eta} \mathcal{Z} f(x+\xi,\omega+\eta).
\end{equation}
It is a well-known property of functions in $W_0(\R)$ that the Zak transform has at least one zero in $[0,1)^2$ and that
\begin{itemize}
    \item if $f \in W_0(\R)$ is even, then
    \begin{equation}\label{eq:zero_even}
        \mathcal{Z} f (x,\omega) = 0, \quad (x,\omega) \in \left(\Z+\tfrac{1}{2}\right)^2,
    \end{equation}

    \item if $f \in W_0(\R)$ is odd, then
    \begin{equation}\label{eq:zero_odd}
        \mathcal{Z} f (x,\omega) = 0, \quad (x,\omega) \in \left(\tfrac{1}{2} \Z \right)^2 \setminus \left(\Z+\tfrac{1}{2}\right)^2.
    \end{equation}
\end{itemize}
This follows from the symmetries of $f$ and the quasi-periodicity of $\mathcal{Z} f$ (see, e.g, \cite[Cor.~4.2]{HorLemVid25}, \cite[\S~5]{Jan88}).
The Poisson summation formula (see, e.g., \cite[Chap.~1.4]{Gro01}) yields, for $f, \widehat{f} \in W_0(\R)$, that the Zak transform satisfies
\begin{equation}\label{eq:Zak_PSF}
     \mathcal{Z}f(x,\omega)
     = e^{-2 \pi i x \omega} \mathcal{Z} \widehat{f}(-\omega,x).
\end{equation}

\section{The Zak transform and Hermite functions}
Using $\F h_n = (-i)^n h_n$, it readily follows from \eqref{eq:Zak_PSF} that
\begin{equation}\label{eq:zero_PSF}
    \mathcal{Z}h_n(0,0) = (-i)^n \mathcal{Z}h_n(0,0) = 0, \quad n \notin 4 \N_0 .
\end{equation}
For $h_{4\ell+1}$ and $h_{4\ell+3}$ ($\ell \in \N_0$) this is also a consequence of the parity, as mentioned above. The fact that $\mathcal{Z} h_{4\ell+2}(0,0) = 0$, was already observed in a different context by Boon, Zak, and Zucker \cite{BooZakZuc83}. For $h_{4\ell}$ we do in general not have a zero at the origin. For example the Zak transform of the Gaussian function $h_0$ only has one simple zero at $(1/2,1/2)$.
\begin{figure}[ht]
    \centering
    \includegraphics[width=0.9\linewidth]{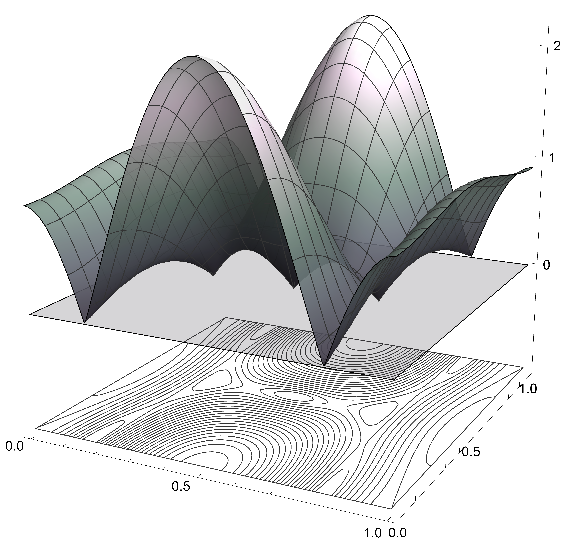}
    \caption{Surface plot of $|\mathcal{Z} \mathcal{D}_{\sqrt{2}}^{-1} h_2(x,\omega)|$ with contour lines on the bottom. $\mathcal{Z} \mathcal{D}_{\sqrt{2}}^{-1} h_2(x,\omega)$ has zeros located at $(x,\omega) \in \{(\frac{1}{4}, \frac{1}{2}), (\frac{1}{2}, \frac{1}{2}), (\frac{3}{4}, \frac{1}{2})\} + \Z^2$. Note that numerically $\mathcal{Z} \mathcal{D}_{\sqrt{2}}^{-1} h_2 (0,\omega) \neq 0$ and $\mathcal{Z} \mathcal{D}_{\sqrt{2}}^{-1} h_2 (x,0) \neq 0$.}
    \label{fig:Zh2D2}
    \vspace*{-12pt}
\end{figure}
In \cite{Lem16}, Lemvig showed, among other results, that (cf.~Fig.~\ref{fig:Zh2D2})
\begin{equation}\label{eq:Zak_h4m+2_Lem}
    \mathcal{Z} (\mathcal{D}_{\sqrt{2}}^{-1} h_{4\ell+2}) \left( \tfrac{1}{4}, \tfrac{1}{2} \right) =
    \mathcal{Z} (\mathcal{D}_{\sqrt{2}}^{-1} h_{4\ell+2}) \left( \tfrac{3}{4}, \tfrac{1}{2} \right) = 0.
\end{equation}
This fact was already noted by Boon, Zak, and Zucker \cite{BooZakZuc83}.

\begin{figure}[ht]
    \centering
    \includegraphics[width=0.8\linewidth]{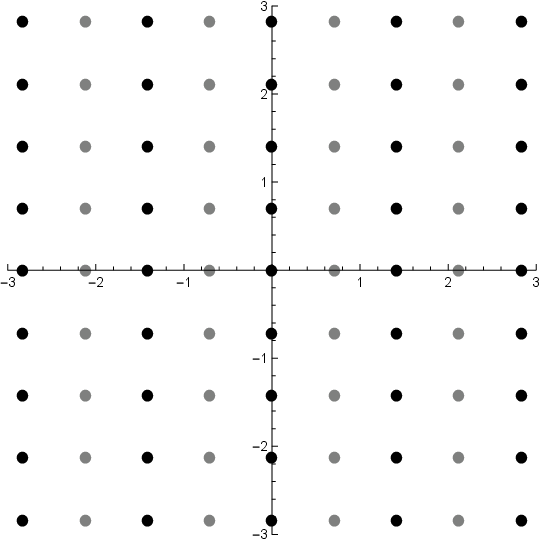}
    \caption{Decomposition of the square lattice $\frac{1}{\sqrt{2}} \Z^2$ into the rectangular lattice $D_{\sqrt{2}} \Z^2$, which has a side ratio of 2, and the shifted copy $D_{\sqrt{2}}(\Z+\frac{1}{2})^2$.}
    \label{fig:lattice}
\end{figure}

To simplify notation for the rest of the section, we set
\begin{equation}
    \widetilde{h}_m = \pi\left(\tfrac{1}{4} + \tfrac{m-1}{2},\tfrac{1}{2} \right) \mathcal{D}_{\sqrt{2}}^{-1} h_{4\ell+2}, \quad m \in \{1,2\}.
\end{equation}
We note that \eqref{eq:Zak_h4m+2_Lem} implies that the multi-window Gabor system
\begin{equation}
    \mathfrak{G}_2 \left( \widetilde{h}_m, \Z^2 \right) = \bigcup_{m=1}^2 \mathcal{G}\left(\widetilde{h}_m, \Z^2\right)
\end{equation}
is not a frame. This stays true by applying any of the above operator-isomorphism pairs or shifting the system by a common shift $\pi(z)$. We have the following unitary equivalences
\begin{align}
    \mathfrak{G}_2 \left( \widetilde{h}_m, \Z^2 \right)
    & \cong \mathcal{G} \left(\mathcal{D}_{\sqrt{2}}^{-1} h_{4\ell+2}, \Z^2 \cup \left((\Z+\tfrac{1}{2}) \times \Z \right) \right)\\
    & \cong \mathcal{G} \Big( h_{4\ell+2}, D_{\sqrt{2}} \left(\Z^2 \cup \left((\Z+\tfrac{1}{2}) \times \Z \right) \right) \Big)
\end{align}
Decomposing the square lattice of density 2 (see Fig.~\ref{fig:lattice}) as
\begin{equation*}
    \tfrac{1}{\sqrt{2}} \Z^2 =
    \left(\sqrt{2} \Z \times \tfrac{1}{\sqrt{2}} \Z \right) \cup
    \left(\sqrt{2} (\Z+\tfrac{1}{2}) \times \tfrac{1}{\sqrt{2}} \Z \right),
\end{equation*}
it follows that $\mathcal{G} \left(h_{4m+2}, \frac{1}{\sqrt{2}} \Z^2 \right)$ is not a frame, which is part of \cite[Thm.~7]{Lem16}. Before proceeding, we state the following lemma (see, e.g., \cite{OzaKatMen99}), which leads to the next result.
\begin{lemma}
    The Hermite function $h_n$ is an eigenfunction of the fractional Fourier transform $\F_r$ with eigenvalue $e^{- i n r}$:
    \begin{equation}
        \F_r h_n = e^{- i n r} h_n, \quad n \in \N_0.
    \end{equation}
\end{lemma}
\begin{corollary}
   $\mathcal{G}\left(h_{4\ell+2}, \tfrac{1}{\sqrt{2}} R_r \Z^2\right)$ is not a frame.
\end{corollary}

\subsection{An alternative proof of Lemvig's result}
\begin{figure}[ht]
    \centering
    \includegraphics[width=0.9\linewidth]{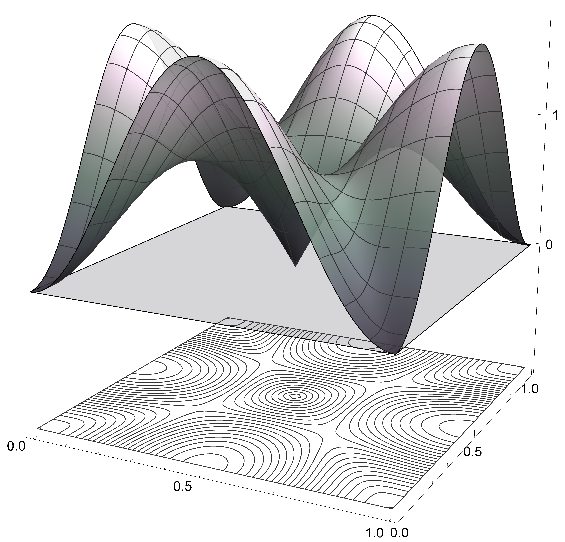}
    \caption{Surface plot of $|\mathcal{Z} h_2(x,\omega)|$ with contour lines on the bottom. Zeros of $\mathcal{Z} h_2(x,\omega)$ are located at $(x,\omega) \in \Z^2 \cup \left(\Z + \frac{1}{2}\right)^2$.}
    \label{fig:Zh2}
\end{figure}

We make use the zeros given by \eqref{eq:zero_even} and \eqref{eq:zero_PSF} (cf.~Fig.~\ref{fig:Zh2}):
\begin{equation}\label{eq:zeros_h4l+2}
    \mathcal{Z}h_{4\ell+2}(0,0) = \mathcal{Z}h_{4\ell+2}(\tfrac{1}{2},\tfrac{1}{2}) = 0
\end{equation}

Note that $\mathcal{Z} (\mathcal{D}_a h_{4\ell+2})(0,0)$ is not necessarily zero and certainly not in general, as indicated by Fig.~\ref{fig:Zh2D2}. Using \eqref{eq:zeros_h4l+2} we easily conclude the following result.
\begin{proposition}
    The Gabor system $\mathcal{G}(h_{4m+2}, \Z^2 \cup (\Z+\frac{1}{2})^2)$ is not a Gabor frame for $L^2(\R)$.
\end{proposition}
This result is just Lemvig's result in disguise. We note that
\begin{equation}
    \Z^2 \cup (\Z+\tfrac{1}{2})^2 = 
    \begin{pmatrix}
        \frac{1}{2} & \frac{1}{2} \\
        -\frac{1}{2} & \frac{1}{2}
    \end{pmatrix} \Z^2 = \tfrac{1}{\sqrt{2}} \mathcal{R}_{\frac{\pi}{4}} \Z^2.
\end{equation}
Thus, this is the same as saying that the Gabor system $\mathcal{G}(\F_{-\frac{\pi}{4}} h_{4m+2}, \frac{1}{\sqrt{2}} \Z^2)$ is not a frame. Since we can safely ignore the phase factor of the window function for the frame property and since $\F_{-\frac{\pi}{4}} h_{4m+2} = e^{\frac{\pi}{4} i (4m+2)} h_{4m+2}$ it follows that we recovered Lemvig's result.

Assume that we do not know Lemvig's result on the zeros of $\mathcal{Z} (\mathcal{D}_{\sqrt{2}}^{-1} h_{4\ell+2})$. By our arguments, however, we know that
\begin{equation}
    \mathcal{G}\left(h_{4\ell+2}, \tfrac{1}{\sqrt{2}} \Z^2\right) \cong \mathcal{G} \left(\mathcal{D}_{\sqrt{2}}^{-1} h_{4\ell+2}, \Z^2 \cup \left((\Z+\tfrac{1}{2}) \times \Z \right) \right).
\end{equation}
Both are not frames. The latter not being a frame implies that
\begin{equation}
    \mathcal{Z}(\mathcal{D}_{\sqrt{2}}^{-1} h_{4m+2})(x_0,\omega_0) = \mathcal{Z}(\mathcal{D}_{\sqrt{2}}^{-1} h_{4m+2})(x_0+\tfrac{1}{2},\omega_0) = 0,
\end{equation}
for some $(x_0,\omega_0) \in [0,1)^2$. Admittedly, we do not have any information on the location of the found zeros, but we proved their existence. We also refer to Janssen's results, in particular on the zeros of the Zak transform \cite[\S~5]{Jan88}, at this point.

\subsection{Another proof that \texorpdfstring{$\mathcal{G}\left(h_2, \tfrac{1}{\sqrt{2}} \Z^2\right)$}{G (h2, sqrt(1/2) Z 2)} is not a frame}
For the special case of the second Hermite function
\begin{equation}
    h_2(t) = 2^{-1/4} (-1 + 4 \pi t^2) e^{-\pi t^2} 
\end{equation}
we will now provide an alternative proof that the Gabor system $\mathcal{G}\left(h_2, \tfrac{1}{\sqrt{2}} \Z^2\right)$ is not a frame. We will again use the Zak transform and \eqref{eq:zero_even}, which says that $\mathcal{Z}f(\frac{1}{2},\frac{1}{2})=0$ when $f$ is an even function. We will then prove the existence of another zero at $(0,0)$ by using the Jacobi identity for theta functions.

We start with writing down the Zak transform
\begin{equation}
    2^{1/4} \mathcal{Z} h_2 (x,\omega) = \sum_{k \in \Z} (-1+4 \pi (k-x)^2) e^{-\pi (k-x)^2} e^{2 \pi i k \omega}.
\end{equation}
Next, we introduce the Jacobi theta-3 function as follows
\begin{equation}
    \vartheta_3(z;\tau) = \sum_{k \in \Z} e^{\pi i \tau k^2} e^{2 \pi i k z}, \quad z \in \mathbb{C}, \, \tau \in \mathbb{H}.
\end{equation}
We refer to \cite{SteSha_Complex_03} or \cite{WhiWat69} for more information. Restricting the arguments, we set
\begin{equation}
    \theta_3(\alpha) = \vartheta_3(0; i \alpha), \quad \alpha > 0.
\end{equation}
The Jacobi identity implies the following functional equation, which also follows from the Poisson summation formula:
\begin{equation}\label{eq:Jacobi}
    \sqrt{\alpha} \, \theta_3(\alpha) = \theta_3\left( \tfrac{1}{\alpha} \right).
\end{equation}
Next, we note that we can write
\begin{equation}
    2^{1/4} \mathcal{Z} h_2 (0,0) = - \theta_3(1) - 4 \theta_3'(1). 
\end{equation}
We differentiate both sides of \eqref{eq:Jacobi} logarithmically to obtain
\begin{align}
    & \quad \frac{1}{2 a} + \frac{\theta_3'(\alpha)}{\theta_3(\alpha)} = - \frac{1}{\alpha^2} \frac{\theta_3'\left( \tfrac{1}{\alpha} \right)}{\theta_3\left( \tfrac{1}{\alpha} \right)}\\
    \Longleftrightarrow & \quad
    \alpha \frac{\theta_3'(\alpha)}{\theta_3(\alpha)} + \frac{1}{\alpha} \frac{\theta_3'\left( \tfrac{1}{\alpha} \right)}{\theta_3\left( \tfrac{1}{\alpha} \right)} = -\frac{1}{2}.
\end{align}
Now, setting $\alpha = 1$, we see that
\begin{equation}
    \theta_3(1) + 4 \theta_3'(1) = 0
    \quad \Longrightarrow \quad
    \mathcal{Z} h_2(0,0) = 0.
\end{equation}
Thus, we found the desired zeros and again conclude that $\mathcal{G}\left(h_2, \frac{1}{\sqrt{2}} R_r \Z^2\right)$ is not a Gabor frame for $L^2(\R)$.

\subsection{Gabor frames: an outlook}
We finish with positive observations, which however need mathematical verification. From Fig.~\ref{fig:Zh2} it is rather obvious that
\begin{equation}
    \mathcal{Z} h_2(x,\omega) \neq 0, \quad (x,\omega) \notin \Z^2 \cup (\Z+\tfrac{1}{2})^2.
\end{equation}
Thus, the double over-sampled Gabor system
\begin{equation}
    \mathcal{G}(h_2, \Z^2 \cup (\Z^2+z)), \quad z \notin (\Z+\tfrac{1}{2})^2
\end{equation}
is potentially a Gabor frame. Moreover, with the presented techniques, it should be possible to find a wealth of Gabor systems $\mathcal{G}(h_2,\L \cup (\L+z))$, with $\L = R_r V_q D_a \Z^2$, which are frames with double over-sampling. This observation is not specific to $h_2$ but is rather valid for any nice $g$ when $\mathcal{Z}g$ only has finitely many zeros. Note that Gröchenig and Lyubarskii~\cite{GroechenigLyubarskii_Hermite_2007} showed that an over-sampling rate of more than $n+1$ suffices for $h_n$ to give a frame over a lattice.


\end{document}